\documentclass[a4paper,12pt]{article}

\usepackage{amsfonts,amsmath,amscd,hyperref,amsthm,makeidx}
\usepackage[dvips]{graphicx}

\newtheorem{thm}{Theorem}
\newtheorem{prop}{Proposition}

\newtheorem{lem}{Lemma}

\newtheorem{rem}{Remark}

\newcommand{\RR}{{\mathbb R}}

\newcommand{\SSS}{{\mathbb  S}}

\newcommand{\abs}[1]{\lvert#1\rvert}
\newcommand{\norm}[1]{\lVert#1\rVert}

\newcommand{\UU}{{\mathcal  U}}

\newcommand{\MM}{{\mathcal  M}}

\newcommand{\dotex}{{\frac{d}{dt}}}

\newcommand{\vv}{{\vec{ v} }}

\begin{document}

\title{A Simple Intrinsic Reduced-Observer for Geodesic Flow}
\author{Silv\`ere Bonnabel
\thanks{S. Bonnabel is with the CAOR, Mines ParisTech. ~{\tt\small silvere.bonnabel@mines-paristech.fr}}}
\date{November 2009}
\maketitle
\begin{abstract}
Aghannan and Rouchon proposed a new design method of asymptotic
observers for a class of nonlinear mechanical systems: Lagrangian
systems with configuration (position) measurements. The (position and velocity) observer is
based on the Riemannian structure of the configuration manifold endowed with
the kinetic energy metric and is intrinsic. They proved local
convergence. When the system is conservative, we propose an intrinsic  reduced order (velocity) observer based on the
Jacobi metric, which can be initialized such that  it  converges exponentially  for \emph{any} initial true velocity. For non-conservative systems the observer can be used
as a complement to the one of Aghannan and Rouchon. More generally
the reduced observer provides velocity estimation for geodesic flow
with position measurements. Thus it can be (formally) used as a fluid flow soft sensor in the case of a perfect incompressible fluid. When the curvature is negative in all planes
the geodesic flow is sensitive to
initial conditions. Surprisingly in this case we have  global exponential convergence and the more unstable the flow is, faster is the
convergence.
\end{abstract}

There is no general method to design asymptotic observers for
observable non-linear systems. Indeed only some
  specific types of
 linearities have been tackled in the literature. In particular
   over the last few years some  work has been devoted  to
  observer design for systems possessing symmetries. \cite{aghannan-rouchon-cdc02,arxiv-07,mahony-et-al-IEEE,bonnabel-mirrahimi-rouchon:07} consider a
  finite-dimensional group of symmetries acting on the state space,
  and \cite{arxiv-08} a left-invariant dynamics on a Lie group. Symmetries generally correspond to invariance
   to some changes of units and frame. Invariance to any change of
   coordinates was raised by \cite{aghannan-rouchon-ieee03} who designed an
intrinsic observer for a class of non-linear systems: Lagrangian
systems  with position (configuration) measurements. The aim is to
estimate the velocity, independently from any nontrivial choice of
coordinates, and of course never differentiate the (noisy) output.
The observer
 was adapted to the specific case of a
left-invariant system on a  Lie
group by~\cite{maithripala2005}. Observer
\cite{aghannan-rouchon-ieee03} is based on the Riemannian structure of
 the configuration manifold endowed with the kinetic energy metric. This geometry had already been used
  in control theory of mechanical systems (see e.g. \cite{lewis-murray00,bullo-lewis-article00}). The
  convergence of the observer is local.

According to the Maupertuis principle, the motion of a conservative
 Lagrangian system is a geodesic flow (motion along a geodesic with constant speed) for the Jacobi metric, intrinsically
 defined using the kinetic and potential energies, up to a time
  reparametrization. In this paper we consider the general problem
  of building a
reduced order velocity observer for geodesic flow on a
Riemannian manifold  with position measurements. A reduced observer is meant to estimate only the unmeasured part of the system's state (here the velocity).
Under some basic assumptions relative to the injectivity radius
(also formulated in \cite{aghannan-rouchon-ieee03}) we
have the following results (Theorem 1). If there is an upper bound $A>0$ on the
sectional  curvature in all planes, choosing $\hat v(0)=0$, the reduced velocity observer always
converges exponentially  to the true velocity, as long as the gain
is larger than  a linear function of $\sqrt A$. Unfortunately the higher the gain is the most sensitive to noise
the observer is. An even better situation occurs when the sectional
curvature is non-positive in all planes: the reduced observer is
globally exponentially convergent for all positive gain.
In fact,   the more negative the curvature is the faster the observer
converges. This feature is surprising enough as negative
curvature implies exponential divergence between two nearby
geodesics, and thus ``amplifies" initial errors.  This  is a major difference with
\cite{aghannan-rouchon-ieee03} who used additional terms precisely
to cancel the effects of (negative) curvature.

For mechanical Lagrangian systems  the observer of Aghannan and Rouchon is only locally convergent. In the absence of external forces the reduced observer provides  an  alternative observer which allows to always estimate the true velocity. When there are external forces, the reduced-observer can be used as a  complement to \cite{aghannan-rouchon-ieee03}. The gain must be chosen  large enough, so that the reduced observer converges before the energy varies significantly. If so, it provides an estimated velocity close to the true one, with which  the observer \cite{aghannan-rouchon-ieee03} can be initialized.

The reduced observer is also applied, formally, to a basic
velocimetry problem: compute the velocity of a perfect
incompressible fluid observing the fluid particles.  The principle
of least action implies that the motion of an incompressible fluid can
be viewed as a geometric flow.  We consider the case of  a
two-dimensional fluid. As the convergence properties of the observer
depend on the sign of the curvature, we will use results and
heuristics of Arnol'd \cite{arnold-book-3}. Following them,  we show that global convergence could be expected  for a large
class of trajectories, since the curvature is positive only in a few
sections. This latter fact also implies instability of the flow, and
Arnol'd  interpretes the difficulty of
weather's prediction as a consequence of this result. Note that the problem tackled is nontrivial,
as the system is  nonlinear, infinite dimensional, and possibly
sensitive to initial conditions.

In Section I we give the general motivations introducing Lagrangian systems on manifolds and Maupertuis' principle. In Section II we introduce the observer. In Section III we consider applications to some mechanical and hydrodynamical systems. In Section IV the convergence in the case of positive constant curvature is illustrated by simulations on the sphere.

\section{Lagrangian systems on manifolds}
Consider the classical mechanical system with $n$ degrees of freedom
 described
by the Lagrangian
$$
\mathcal   {L}(q,\dot q) = \frac{1}{2} g_{ij}(q) \dot{q}^i \dot{q}^j
- U(q)
$$
where the generalized positions $q\in \mathcal M$ are written in the local
coordinates $(q^i)_{i=1...n}$,  $g(q) = (g_{ij}(q))_{i=1\ldots
n,j =1\ldots n}$ is a Riemaniann metric on the configuration space
$\mathcal M$, and $U: \mathcal M \mapsto \RR$ is the potential energy. The
Euler-Lagrange equations write in the local coordinates
\begin{align}\label{Euler_lagrange:eq}
{\frac{d}{dt}} \left(\frac{\partial}{\partial\dot{q}^i}\mathcal
{L}\right) = \frac{\partial}{\partial q^i}\mathcal   {L}, \quad
i=1,...,n
\end{align}
One can prove using $\frac{\partial
g^{ik}}{\partial q^l}g_{jk}=-g^{ik}\frac{\partial g_{jk}}{\partial
q^l}$ where $g^{il}$  are components of $g^{-1}$ that
\eqref{Euler_lagrange:eq} writes
\begin{align} \label{dyn:eq}
  \begin{array}{rl}
\ddot{q}^i &= - \Gamma^i_{jk}(q) \dot{q}^j \dot{q}^k +\frac{\partial
}{\partial q^i}U
  \end{array}
\end{align}
where the Christoffel symbols  $\Gamma^i_{jk}$ are given by $
\Gamma^i_{jk} = \frac{1}{2} g^{il}\left(\frac{\partial g_{lk}
}{\partial q^j}+\frac{\partial g_{jl} }{\partial q^k} -
\frac{\partial g_{jk} }{\partial q^l}\right) $ (see e.g. \cite{abraham-marsden-book}) . A curve $\gamma(t)$
which is a critical point of the action $$S(\gamma)=\int_0^T
L(\gamma(t),\dot\gamma(t) dt)$$ among all curves with fixed
endpoints satisfies the Euler-Lagrange equations
\eqref{Euler_lagrange:eq}.
\subsection{Lagrangian system in a potential field}\label{maupertui:sec}
Consider a conservative Lagrangian system evolving in an admissible
region $\{q\in \mathcal M: U(q)<E\}$. The energy of the system
$E=T(q,\dot q)+U(q)=\frac{1}{2} g_{ij}(q) \dot{q}^i \dot{q}^j+U(q)$
is fixed. According to the Maupertuis principle of least action (see
e.g. \cite{arnold-book-3}),  in the Riemannian geometry defined by the Jacobi metric $\hat
g^{ij}(q)=2(E-U(q))g^{ij}(q)$ and the natural parameter $\tau$ such
that $\frac{d\tau}{dt}=2(E-U(q(t)))$, the geodesic flow is a solution of the equation
of motion \eqref{dyn:eq}. Indeed if the $\hat\Gamma^i_{jk}$ are the
Christoffel symbols associated to the metric $\hat g$ we have
\begin{align}\label{}
\frac{d^2}{d\tau^2}{q}^i + \hat\Gamma^i_{jk}(q) \frac{d}{d\tau}{q}^j
\frac{d}{d\tau}{q}^k =0
\end{align}which writes intrinsically $\hat\nabla_\frac{dq}{d\tau}\frac{dq}{d\tau}=0$  and defines the geodesic flow ($\hat\nabla$ is the Levi-Civita covariant differentiation of the Jacobi metric).
\subsection{Geodesic flow and holonomic constraints}\label{holonomic:sec}
A material particle constrained to lie on a manifold moves along a
geodesic \cite{arnold-book-3}. Indeed
 $E=T$, $\hat g=2Eg$,  $d\tau=2Edt$  ensure the energy $T$ is fixed. According to Maupertuis'
principle the motion
  minimizes
 $\int_\gamma\sqrt{\hat g_{ij}\frac{d}{d\tau}{q}^i
\frac{d}{d\tau}{q}^j}d\tau=(1/{\sqrt{2E}})\int_\gamma\sqrt{
g_{ij}\frac{d}{dt}{q}^i \frac{d}{dt}{q}^j}dt$ which is proportional
to the geodesic length in the metric $g$. More generally an inertial motion of a Lagrangian system with $k$ holonomic constraints can be viewed as the inertial motion of a particle constrained to lie on a submanifold of dimension $n-k$ (see e.g. \cite{arnold-book-3} p 90). A conservative Lagrangian system in a potential field with holonomic constraints satisfies the Maupertuis' principle on the configuration submanifold of dimension $n-k$.

\section{An intrinsic reduced-observer}
Let us build an observer to estimate the velocity $\dot q$ of a point $q$  moving  along the geodesics of $\MM$ with constant speed, when the position $q$ is measured (with noise). First suppose  $\mathcal M=\mathbb{R}^n$ endowed with Euclidian metric.  Let $\dot q=v$ and $\dot v=0$. For such a linear system a Luenberger reduced dimension  observer with arbitrary dynamics can be constructed \cite{Luenberger}. The goal is to estimate only the part of the state that is not directly measured. An auxiliary variable $\xi$, which is a combination between the unmeasured part of the state and the output, is generally introduced: \begin{align}\xi=q-\lambda v\qquad(\text{and~thus}~~\dot\xi=v)
\end{align}To estimate $\xi$ and $\dot\xi$ consider the reduced observer:
\begin{align}\label{plan}
\dotex \hat\xi=-\frac{\hat\xi-q}{\lambda}
\end{align}
It can be interpreted as a simple pursuit algorithm with proportional feedback. Let $\hat v=\dotex\hat\xi$.
Let us prove $\hat\xi-\xi\rightarrow0$ and $\hat v - v\rightarrow 0$.  We have $\dotex {\hat
v}=-(1/\lambda)(\hat v - v)$ implying  $\hat
v\rightarrow v$ for $\lambda>0$. As $\hat\xi=q-\lambda\hat v$, we have    $\hat\xi-\xi\rightarrow0$, and $\hat\xi$ is asymptotically  moving behind $q$ at \emph{fixed} distance $\lambda \norm{\dot q}$. If $\mathcal M$ is any Riemannian manifold consider \begin{align}\label{obsreduit}
\dotex \hat\xi=-\frac{1}{2\lambda} ~\overrightarrow {grad}_{\hat\xi}
~D^2(\hat\xi,q),\quad \lambda>0
\end{align}
where $D (\hat \xi,q)$ is the geodesic distance between $\hat \xi$ and $q$. If $D(\hat\xi,q)$ is smaller than the  injectivity radius at $q$, then \eqref{obsreduit} means that $\dotex \hat\xi$ is a vector  which is tangent to the geodesic linking $\hat \xi$ and $q$, and whose norm is proportional to $D(\hat\xi,q)$. The dynamic does not depend on any choice of local coordinates in $\RR^n$,
and is a generalization of \eqref{plan}. We want to prove that $D(\hat \xi,\xi)\rightarrow 0$ where $\xi$ is a point following $q$ at distance $\lambda\norm{\dot q}$ on the geodesic $\{q(t):t>0\}$. The parallel transport $\mathcal{T}_{//\hat \xi\rightarrow q}$ of $\dotex
\hat \xi$ to the tangent space at $q$ along the geodesic joining $\hat \xi$ and $q$  is an estimation of
$v=\dot q$.
\begin{align}\label{obs:vitesse:eq}
\hat v=\mathcal{T}_{//\hat \xi\rightarrow q}\dotex \hat \xi \end{align}
\begin{thm}Let $\MM$ be a Riemannian manifold. Let $T<\infty$. Let $t\mapsto q(t)\in\MM$ satisfy $\nabla_{\dot q} \dot q=0$ for $t\in[0,T]$.  Let $\xi(t)=\exp_{q(t)}(-\lambda {\dot q})$. Consider the observer \eqref{obsreduit}. Let \eqref{conv:ineq} be the inequality\begin{align}\label{conv:ineq}D(\hat \xi(t), \xi(t))\leq e^{-\frac{1}{\lambda}t} D(\hat \xi(0), \xi(0)) \quad \forall t\in[0,T]\end{align}
\begin{itemize}
\item Suppose the Riemannian curvature is non-positive in all planes. If  for all $t\in[0,T]$, $D(\hat\xi(t),q(t))$ is bounded by the injectivity radius $I(t)$ at $q(t)$ (i.e. there exists a unique geodesic joining $\hat\xi$ and $q(t)$), \eqref{conv:ineq} is true for all $\lambda>0$. When the manifold is complete and simply-connected (Hadamard manifold), the   injectivity radius is infinite (Cartan-Hadamard theorem) and \eqref{conv:ineq} is always true. In particular $\hat \xi(0)$ can be chosen arbitrarily. Moreover for all $t\in[0,T]$ \begin{align}\label{speed:conv:eq1}\lambda\norm{\hat v(t)- \dot q(t)}\leq D(\hat\xi(t),\xi(t))\leq e^{-\frac{1}{\lambda}t} D(\hat \xi(0), \xi(0)) \end{align}
\item
Suppose the sectional curvature is bounded from above by $A>0$.  \eqref{conv:ineq} is true as long as the distance $D(\hat \xi(t),q(t))$ remains bounded by $\max(\frac{\pi}{4\sqrt A},I(t))$  for all $t\in[0,T]$. If the manifold is simply connected and  ${\lambda}>\frac{\pi}{4\norm{\dot q}_g\sqrt A}$, \eqref{conv:ineq} is true as soon as $D(\hat \xi(0),q(0))<\frac{\pi}{4\sqrt A}.$  Moreover in this case we have exponential convergence in polar coordinates for all $t\in[0,T]$: \begin{equation}\begin{aligned}\label{speed:conv:eq2} &\lambda|~\norm{\hat v(t)}-\norm{\dot q}~|\leq D(\hat \xi(t),\xi (t)),\\&0\leq\sin(\alpha_A(t))\leq \frac{ \sqrt A~ }{\sin(\sqrt A~  \lambda\norm{\dot q})}D(\hat\xi(t),\xi(t))\end{aligned}\end{equation} where  $\alpha_A(t)\rightarrow 0$ and the   angle $\alpha(t)$  between $\hat v(t)$ and $\dot q(t)$  satisfies $0\leq\alpha(t)\leq\alpha_A(t)\leq\pi$.
\end{itemize}
\end{thm}
The convergence \eqref{conv:ineq} of the observer's state variable $\hat\xi$ is not sufficient to prove that $\hat v$ converges. Indeed the estimated velocity $\hat v$ is linked to $\hat\xi$ via a non-linear geometric transformation. Yet geometry of  triangles on curved surfaces will allow to prove \eqref{speed:conv:eq1} and \eqref{speed:conv:eq2}.
\begin{proof}
The proof utilizes two differential geometry lemmas. Lemma \ref{kupka:lem} is a consequence of Synge's lemma (see e.g. \cite{stoker-book} p 316) for which  a direct demonstration is proposed.
\begin{lem}\label{contraction:lem}
Let $\MM$ be a smooth Riemannian manifold. Let $P\in \MM$ be fixed. On the subspace of $\MM$ defined by the injectivity radius at $P$ we consider
\begin{align}\label{variete:dyn}
\dotex x=-\frac{1}{2\lambda} ~\overrightarrow {grad}_{x}
~D^2(P,x)\quad \lambda>0
\end{align}If the sectional curvature is non-positive in all planes, the dynamics is a contraction in the sense of \cite{slotine-auto}, i.e, if $\delta x$ is a virtual displacement at fixed $t$ we have
\begin{align}\label{cont:eq}\dotex \|\delta x\|_g^2\leq -\frac{2}{\lambda}\|\delta x\|_g^2\end{align} where
$\|~\|_g~$ is the norm associated to the metric $g$.
If the sectional curvature in all planes is upper bounded from above by $A>0$, \eqref{cont:eq} holds for $D(P,x)<\pi/(4\sqrt A)$.
\end{lem}
\begin{proof}
The virtual displacement is defined \cite{slotine-auto} as a linear
tangent differential form, and can be viewed by duality as a vector
of $TM|_x$. Let us define a surface $\Sigma$. Let $\gamma_0$ be the
geodesic joining $P$ to $x$. Consider $x_\epsilon=\exp_{x}(\epsilon~
\delta x)\in\MM$. It is linked to $P$ by a geodesic, say
$\gamma_\epsilon$. Up to second order terms in $\epsilon$ we have
$\dot \gamma_\epsilon(0)-\dot \gamma_0(0)=\epsilon u$ where $u$ is a
tangent vector at $P$. The directions defined by $\gamma_0$ and $u$
at $P$ span a 2-plane tangent at $P$. All the geodesics having a
direction tangent to this 2-plane at $P$ span a smooth surface
$\Sigma$  embedded in $\MM$ which inherits the Riemannian metric
$g$. We have $x\in \Sigma$ and $\delta x\in T_x\Sigma\subset
T_x\MM$. $\Sigma$ is invariant under the flow $\eqref{variete:dyn}$,
as the gradient term is tangent to the geodesics heading towards
$P$. Indeed, let $\gamma$ be parameterized by the arclength
$\sigma$, and let  $\sigma_0=D(P,x)$. The squared distance increases
the most in the direction of the geodesics. Thus the gradient  tangent to the geodesic. We have
$D^2(P,\gamma(\sigma_0+\epsilon))=(\sigma_0+\epsilon)^2$. Up to
second order terms
$D^2(P,\gamma(\sigma_0+\epsilon))=D^2(P,x)+\epsilon\langle\overrightarrow
{grad}_{x} ~D^2(P,x),\frac{d\gamma}{d\sigma} (\sigma_0)\rangle_g$.
The norm of the gradient is thus is $2D(P,x)$.

Following \cite{stoker-book} (p 177) we use specific coordinates on
$\Sigma$  called ``polar coordinates". Let $e_1,e_2$ be an euclidian
frame of $T_P\Sigma$ for the inherited metric and $e_1$ be tangent to
$\gamma_0$. We define $\Phi: (\sigma,\theta)\mapsto \exp_P(\sigma
\cos\theta e_1+\sigma\sin\theta e_2)$. $\Sigma$ is parameterized by
$\sigma$, the geodesic length  to $P$, and $\theta$, the angle  in
$T_P\Sigma$ with $e_1$. In the polar coordinates, the  elementary
length is given by
$$ds^2=d\sigma^2+G(\sigma,\theta)d\theta^2$$ and $G$ satisfies the initial conditions $\sqrt G=0$ and $\frac{\partial \sqrt G}{\partial
\sigma}=1$ at $\sigma=0$. According to a classical result \cite{stoker-book} the Gauss curvature at the point  $u=\Phi(\sigma,\theta)$ is given by
 $K(\sigma,\theta)=\frac{-1}{\sqrt {G(\sigma,\theta)}}\frac{\partial^2
\sqrt {G(\sigma,\theta)}}{\partial \sigma^2}$. We will prove (lemma \ref{kupka:lem}) that the Gaussian curvature at $u=\Phi(\sigma,\theta)\in \Sigma$ is less than the sectional curvature in the tangent plane to $\Sigma$ at $u$: $K(\sigma,\theta)\leq K_{sec}(T_u\Sigma)$.

Suppose $K_{sec}(T_u\Sigma)\leq 0$. It implies $K(\sigma,\theta)\leq 0$. Along $\gamma$ we have
\begin{equation}\begin{aligned}\label{ineq:eq}\frac{\partial^2 G(\sigma,\theta)}{\partial \sigma^2} &=\frac{\partial }{\partial
\sigma}(2\sqrt {G(\sigma,\theta)}\frac{\partial \sqrt {G(\sigma,\theta)}}{\partial \sigma})\\&=2\left((\frac{\partial \sqrt {G(\sigma,\theta)}}{\partial \sigma})^2-G(\sigma,\theta)K(\sigma,\theta)\right)\geq 0\end{aligned}\end{equation} and
thus $\frac{\partial }{\partial \sigma}(\sigma\frac{\partial
G}{\partial \sigma})\geq \frac{\partial G}{\partial \sigma}$ which yields by integration
$\sigma\frac{\partial G}{\partial \sigma}\geq G$ since $G(0,\theta)=0$. In the polar coordinates the dynamics \eqref{variete:dyn} reads
$$\dot \sigma=-\frac{1}{\lambda}\sigma;\qquad \dot\theta=0$$Indeed we already stated that the gradient is tangent to the geodesic, thus $\dot\theta=0$, and \eqref{variete:dyn} becomes a one-dimensional dynamics along the geodesic, and as $D^2(P,x)=\sigma^2$ we have $\norm{\overrightarrow {grad}_{x}
~D^2(x,P)}_g=2\sigma$. Writing $\|\delta x\|^2=\alpha^2\delta\sigma^2+\beta^2G(\sigma,0)\delta\theta^2$ we have along the geodesic $\gamma_0$ (parameterized by $\sigma$ and $\theta=0$) the following inequality, proving \eqref{cont:eq}.
\begin{align}\label{contr:eq}\dotex \|\delta x\|^2=-2\frac{\alpha^2}{\lambda}\delta\sigma^2
-2\beta^2\frac{\sigma}{\lambda}\frac{\partial {G(\sigma,0)}}{\partial
\sigma}\delta\theta^2\leq -\frac{2}{\lambda}\|\delta x\|^2\end{align}

Suppose now $K_{sec}(T_u\Sigma)\leq A$. Let $z(\sigma)=\sqrt {G(\sigma,0)}$. We have $z''=-K(\sigma,0)z$,  $z(0)=0,~z'(0)=1$, with $K(\sigma,0)\leq A$.  The Sturm comparison theorem allows to compare $z$ to the solution of equation $y''=-Ay$,  $y(0)=0,~y'(0)=1$, i.e. $y(\sigma)=\sin(\sqrt A \sigma)/\sqrt A$. A Taylor expansion in $0$ shows  there exists $\mu>0$ such that $z(\sigma)/z'(\sigma)<  y(\sigma)/y'(\sigma)$ for $0<\sigma\leq \mu$. It is proved in \cite{arnold-book-1} (Sturm Comparison theorem) this implies $z'(\sigma)>0$  and $z(\sigma)/z'(\sigma)<  y(\sigma)/y'(\sigma)$ for $0<\sigma < \pi/(2\sqrt A)$. Indeed the last inequality is based on the fact that  $(z/z')'=1+K(z/z')^2\leq 1+A(z/z')^2$ and thus $z/z'$ can never ``overtake" $y/y'$ (see \cite{arnold-book-1}). Thus for $0\leq\sigma\leq\pi/(4\sqrt A)$ we have $z(\sigma)\leq z'(\sigma) \tan(\pi/4)/\sqrt{A}$, and thus $z'(\sigma)^2-Az(\sigma)^2\geq 0$. Thus for $\sigma=D(P,x)\leq\pi/(4\sqrt A)$, \eqref{ineq:eq} is true (with $\theta\equiv 0$) and \eqref{cont:eq} holds.

%
%
%

\end{proof}
\begin{lem}\label{kupka:lem}Let $\MM$ be a smooth manifold. Let $P\in\MM$. Let E be a two-dimensional vectorial space of $T_P\MM$. Let $\omega$ be a neighborhood of $P$ in $E$ such that the restriction of the exponential map $\rho$ to $\omega$ is a diffeomorphism in $\MM$.  $\Sigma=\rho(\omega)$ is submanifold of dimension 2. Its Gaussian curvature at any $u\in \Sigma$ is less than the sectional curvature in the tangent plane to $\Sigma$ at $u$: $$K(u)\leq K_{sec}(T_u\Sigma)$$
\end{lem}
\begin{proof}
The proof is based on computations due to Ivan Kupka. Let $\Sigma$ be the surface of lemma \ref{contraction:lem} and $\Phi$ the map associated to the polar coordinates. The metric inherited by $\Sigma$ writes\begin{align*}ds^2&=\norm{\frac{\partial\Phi}{\partial\sigma}dr+\frac{\partial\Phi}{\partial\theta}d\theta}^2_g=\norm{\frac{\partial\Phi}{\partial\sigma}}_g^2dr^2\cdots\\&\quad+2<\frac{d\Phi}{\partial\sigma},\frac{\partial\Phi}{\partial\theta}>_gd\sigma d\theta+\norm{\frac{\partial\Phi}{\partial\theta}}^2_gd\theta^2\end{align*}For fixed $\theta$ the curve $\gamma_\theta: \sigma\mapsto \Phi(\sigma,\theta)$ is a geodesic and thus $\norm{\frac{\partial\Phi}{\partial\sigma}}^2_g=1$. Let  $J(\sigma,\theta)=\frac{\partial\Phi}{\partial\theta}$. For fixed $\theta$, $J_\theta: \sigma\mapsto J(\sigma,\theta)$ is a Jacobi field along $\gamma_\theta$. Moreover $J(0,\theta)=-\sin\theta e_1+\cos\theta e_2$ is orthogonal to this geodesic at $P$. It is well known (Jacobi field properties) that it implies $J$ is orthogonal to $\gamma_\theta$ at any point. Thus  $\langle\frac{\partial\Phi}{\partial\sigma},\frac{\partial\Phi}{\partial\theta}\rangle_g=0$ and $ds^2=d\sigma^2+\norm{J(\sigma,\theta)}_g^2d\theta^2$, and the Gaussian curvature  is given by $K(\sigma,\theta)=\frac{-1}{\norm{J(\sigma,\theta)}_g}\frac{\partial^2
{\norm{J(\sigma,\theta)}_g}}{\partial \sigma^2}$ (see \cite{stoker-book}). Consider the Levi-Civita  covariant differentiation $\nabla$ of the metric $g$. We have $\frac{\partial
{\norm{J(\sigma,\theta)}_g}}{\partial \sigma}=\frac{\langle\nabla_\sigma J(\sigma,\theta),J(\sigma,\theta)\rangle}{\norm{J(\sigma,\theta)}_g}$ and
\begin{align*}&\frac{\partial^2
{\norm{J(\sigma,\theta)}_g}}{\partial \sigma^2}= \frac{\langle\nabla_\sigma^2 J(\sigma,\theta),J(\sigma,\theta)\rangle}{\norm{J(\sigma,\theta)}_g}\cdots\\&\qquad+\frac{\langle\nabla_\sigma J(\sigma,\theta),\nabla_\sigma J(\sigma,\theta)\rangle}{\norm{J(\sigma,\theta)}_g}-\frac{\langle\nabla_\sigma J(\sigma,\theta),J(\sigma,\theta)\rangle^2}{\norm{J(\sigma,\theta)}_g^3}\end{align*}According to the Jacobi equation we have $\nabla_\sigma^2 J(\sigma,\theta)+R(J(\sigma,\theta),\frac{\partial\Phi(\sigma,\theta)}{\partial\sigma})\frac{\partial\Phi(\sigma,\theta)}{\partial\sigma}=0$. Thus the Gaussian curvature of $\Sigma$ satisfies \begin{align*}&K(\sigma,\theta)=K_{sec}(T_u\Sigma)\cdots\\&\qquad+\frac{\langle\nabla_\sigma J(\sigma,\theta),J(\sigma,\theta)\rangle^2-\norm{J(\sigma,\theta)}_g^2{\norm{\nabla_\sigma J(\sigma,\theta)}_g^2}}{\norm{J(\sigma,\theta)}_g^4}\end{align*} where $K_{sec}(T_u\Sigma)=\frac{<R(J(\sigma,\theta),\frac{\partial\Phi(\sigma,\theta)}{\partial\sigma})\frac{\partial\Phi(\sigma,\theta)}{\partial\sigma},J(\sigma,\theta)>}{\norm{J(\sigma,\theta)}^2_g}$ is the value of the sectional curvature on the tangent plane to $\Sigma$ at $u$, and where we used that $J(\sigma,\theta)$ is orthogonal to $\frac{\partial\Phi(\sigma,\theta)}{\partial\sigma}$. Cauchy-Schwarz implies that the fraction above is negative and $K(\sigma,\theta)\leq K_{sec}(T_u\Sigma)$.
\end{proof}

\begin{figure}[htb]
\centering
\includegraphics*[width=1\textwidth]{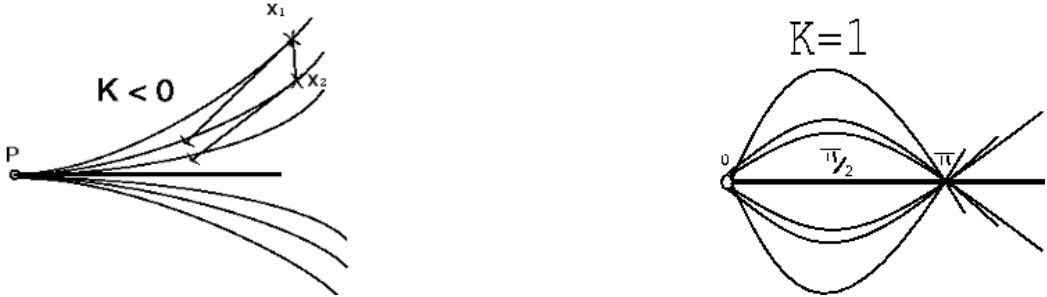}
  \caption{Left: Geodesic deviation on a manifold of negative curvature.  \eqref{variete:dyn} writes in polar coordinates $(\sigma,\theta)$: $\dot \sigma=-\frac{\sigma}{\lambda}; \dot\theta=0$. The distance $\norm{\delta x}$ between neighbors $x_1$ and $x_2$ decreases at a rate at least $\frac{1}{\lambda} $. Right: Geodesic deviation on the sphere.}
    \label{deviation:lem:fig}
\end{figure}
Suppose that for all $t\in[0,T]$ $D(\hat\xi,q(t))$ is bounded by the injectivity radius, as well as by $\frac{\pi}{4\sqrt A}$ in case of positive sectional curvature.  At each time $t$ letting $P=q(t)$ we see that \eqref{obsreduit} is the same as   \eqref{variete:dyn}. Lemma \ref{contraction:lem} thus  proves that  we have the property \eqref{cont:eq} at time $t$, with a contraction rate $\frac{2}{\lambda}$ independent from $t$. Thus Lemma \ref{contraction:lem} used at every $t$  proves that \eqref{obsreduit} is a contraction as defined in  \cite{slotine-auto}, \cite{aghannan-rouchon-ieee03}.  Using
the contraction  interpretation  in the appendix
of \cite{aghannan-rouchon-ieee03} we see that if $\hat \xi_1,\hat
\xi_2$ are solutions of \eqref{obsreduit}  we have
$$D(\hat \xi_1(t),\hat \xi_2(t))\leq e^{-\frac{1}{\lambda}t} D(\hat
\xi_1(0),\hat \xi_2(0)) \quad \forall t\in[0,T]$$The system ``forgets" its
initial condition. So \eqref{conv:ineq} holds if $\xi(t)$ is
a solution of  \eqref{obsreduit}. This is true since $0=\dotex
D(\xi,q)=\norm{v}_g-\norm{\dot\xi}_g=\norm{v}_g-\frac{1}{\lambda}D(\xi,q)$.

Under the basic assumption that $D(\hat\xi,q(t))\leq I(t)$, we have just proved that when the sectional curvature is nonpositive in all planes, \eqref{conv:ineq} is  true for any initial condition. When   the sectional curvature is bounded from above by $A$, we proved at  Lemma \ref{contraction:lem} that \eqref{conv:ineq} holds if for all $t>0$ $D(\hat\xi,q(t))\leq\frac{\pi}{4\sqrt A}$, i.e. $\hat\xi$ remains in the contraction region. Thus, the bound on $\lambda$ is meant to make the contraction region a  trapping region. Indeed $\dotex D(\hat \xi,q)\leq\langle\text{grad}_qD(\xi,q),v\rangle-\norm{\dotex\hat\xi}_g\leq\norm{v}_g-\frac{1}{\lambda}D(\hat \xi,q)$. Thus $D(\hat \xi,q)=\frac{\pi}{4\sqrt A}$ implies $\dotex D(\hat \xi,q)<0$ if ${\lambda}>\frac{\pi}{4\norm{v}_g\sqrt A}$. Thus for ${\lambda}>\frac{\pi}{4\norm{v}_g\sqrt A}$ the vector field is pointing inside the contraction region. In particular if $D(\hat \xi(0),q(0))<\frac{\pi}{4\sqrt A}$ we have $D(\hat \xi(t),q(t))<\frac{\pi}{4\sqrt A}$ for all $t>0$ and \eqref{conv:ineq} holds.

Now that we have proved the exponential convergence of $D(\hat\xi,\xi)$ we can focus on the convergence of $\hat v$ towards $\dot q$.  $(q,\xi,\hat\xi)$ is a geodetic triangle $T$. Indeed as $q,\xi,\hat\xi$ are assumed to be in a ball of radius $\frac{\pi}{4\sqrt A}$, $T$ is well-defined and the angles are less than $\pi$. The length of the sides are: $D(q,\xi)=\lambda \norm{\dot q}$, which is fixed,  $D(q,\hat\xi)=\lambda\norm{\hat v}$, and $D(\hat\xi(t),\xi(t))\leq\exp(-t/\lambda)D(\hat\xi(0),\xi(0))$. In the Euclidian case $(K\equiv 0)$, there is an homothety between $\hat v,\dot q$ and the sides of the triangle and we have $\lambda\norm{\hat v-\dot q}=D(\hat\xi,\xi)$, proving \eqref{speed:conv:eq1}. As $T$ and its sides belong to the surface $\Sigma$  defined  before, one can apply Alexandrov's theorem \cite{stoker-book} stating that, if $A$ is an upper bound on the curvature, the angle $\alpha$ between $\dot q$ and $\hat v$ (the triangle angle at $q$) is less than the angle $\alpha_A$ corresponding to the case of constant curvature $K\equiv A$. Thus when the curvature is nonpositive in all planes, $\alpha$ is less than in the Euclidian case and \eqref{speed:conv:eq1} is proved. When the curvature is upper bounded by $A>0$, $\alpha$ is less that $\alpha_A$ verifying the spherical law
 of sines: $\sin(\alpha_A)=\sin(\beta)\sin( \sqrt A~ D(\hat\xi,\xi))/\sin(\sqrt A~  D(q,\xi))$, where $\beta$ is the opposite angle to the side linking $q$ and $\xi$ \cite{stoker-book}. To prove \eqref{speed:conv:eq2} we used $0\leq\sin\beta\leq 1$. The first part of \eqref{speed:conv:eq2} is the spherical triangular inequality. $\alpha_A\rightarrow \pi$ is impossible as $D(\hat\xi,\xi)\rightarrow 0$.

\end{proof}

\section{Applications}
\subsection{Lagrangian mechanical system}
 \begin{prop}Consider any Lagrangian system in a potential field in the admissible region defined by $U<E$.  The observer  \begin{align}\label{maupertuis:obs:eq}\frac{d}{dt} \hat \xi=-\frac{1}{\lambda}(E-U(q)) ~\overrightarrow {grad}_{\hat \xi}
~D_{\hat g}^2(\hat \xi,q)\end{align} is such that the Theorem  1 is valid in the Maupertuis time and  Jacobi metric.\end{prop}
\begin{proof}
One can apply the Maupertuis' principle (see section
\ref{maupertui:sec}). In Maupertuis' time
$\tau=\int_0^t2(E-U(q(t))dt$, the motion is a geodesic flow on the
configuration space with modified metric $\hat g$, with
$\norm{v}_{\hat g}=1$ as $\hat
g_{ij}(q)\frac{dq^i}{d\tau}\frac{dq^j}{d\tau}=2(E-U)g_{ij}(q)\frac{dq^i}{dt}\frac{dq^j}{dt}(\frac{dt}{d\tau})^2=1$.
The observer defined by $\frac{d}{d\tau} \hat\xi=-\frac{1}{2\lambda}
\overrightarrow{grad}_{\hat \xi} {D_{\hat g}}^2 (\hat \xi,q)$,
$\lambda>0$  where $D_{\hat g}$ is the distance associated to Jacobi
metric, is such that $\hat v=\mathcal{T}_{//\hat \xi\rightarrow
q}\frac{d}{d\tau} \hat \xi$ is an estimation of $\frac{d}{d\tau}q$.
\end{proof}

For instance, R. Montgomery studied in a recent paper \cite{montgomery} the Newtonian equal-mass three bodies problem, with zero momentum and when the potential is taken equal to $1/r^2$: the Jacobi metric has negative curvature everywhere (except at two points). The reduced observer \eqref{maupertuis:obs:eq} is thus globally convergent for a three bodies system which is sensible to initial conditions.
\begin{rem}For a conservative system, the total energy $E$  needs to be known to compute the Jacobi metric. But no information about the  direction of the velocity is required.
\begin{figure}[htb]
\centering
\includegraphics*[width=1\textwidth]{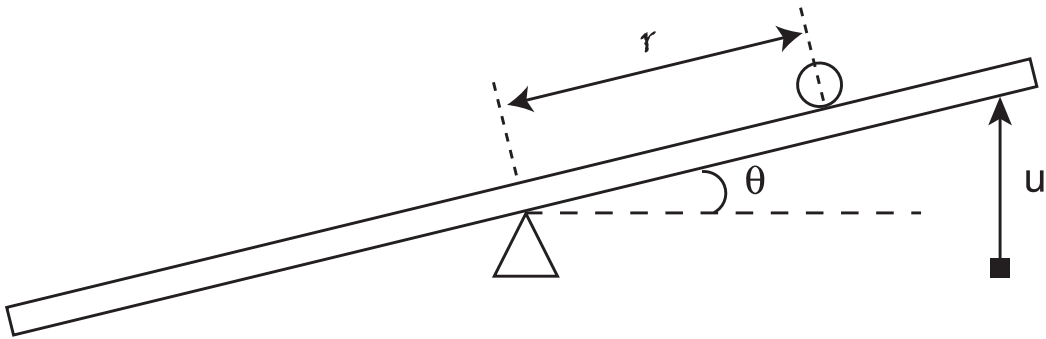}
  \caption{Ball and Beam}
    \label{ballandbeam}
\end{figure}
\end{rem}
\begin{rem}Let us consider now a non-conservative system:  the ball and beam of  \cite{aghannan-rouchon-ieee03} with a torque control $u$  (see fig \ref{ballandbeam}). The observer  \cite{aghannan-rouchon-ieee03} is only {locally} convergent. Observer \eqref{obsreduit} can be used complementarily to provide a globally convergent estimator with the following little experiment. A some time $t_0$ maintain $u\equiv0$ (no control) and set $\hat\xi(t_0)=q(t_0)$. The characteristic time of convergence of the observer \eqref{obsreduit} is $\tau=\lambda$ in the Maupertuis time. After a few  $\tau$ the observer \eqref{obsreduit} provides the observer of Aghannan and Rouchon an initial estimation of the velocity close to the true one and from that moment $u$ can vary freely again: \cite{aghannan-rouchon-ieee03} converges.  As the observer allows to identify the direction of the velocity, it is more interesting to use it for a 3D ball and beam problem in which the beam is replaced with a plate fixed at a point, rotating around two horizontal axis (so that two angles are involved). \end{rem}

\subsection{Motion of a perfect incompressible fluid}The goal of this section is to show that the reduced observer could possibly be applied to more complicated systems. No formal proof is given but only heuristic discussions. The observer could be used in particle velocimetry as a (soft) velocimeter  for a  flow seeded with observable particles and modeled by Euler equations.
\subsubsection{A reduced observer}Let us first introduce some results and notations of \cite{arnold-book-3,arnold-fourier,rouchon-EJM/B91}. Let $\Omega$ be a domain of $\RR^3$ bounded by a surface $\delta\Omega$. Let $\vec{v}$ be the velocity field of an ideal incompressible perfect fluid with density $\rho$ which fills the domain $\Omega$. The motion is described  by the Euler equation\begin{align}\label{euler:eq}
\dotex \vv+(\vv.\nabla)\vv=-\frac{1}{\rho}\nabla p
\end{align}where $p$ is the pressure. Let SDiff~$\Omega$ be the Lie group of all diffeomorphisms that preserve the Euclidian volume. Its Lie algebra  $\UU$ is the set or all vector fields of $\Omega$ of null divergence, and tangent to the boundary $\delta\Omega$.  Consider the scalar product on the Lie algebra\begin{align}\label{kinetic:fluid:eq}\forall \vv,\vec{w}\in\UU,\quad<\vv,\vec{w}>={\rho}\int\int\int_\Omega\vv(x).\vec{w}(x)dx\end{align} Let $\vv(t)\in\UU$ be a solution of \eqref{euler:eq}. Let $\phi_t^\vv(x)$ be the position  at time $t$ of a fluid particle initially at $x$, i.e. obtained by integration on $[0,t]$ of the system $\frac{d}{ds}z=\vv(s,z),~z(0)=x$. $\phi_t^\vv$ is a diffeomorphism for any $t>0$, and the motion of the fluid is described by a curve $t\mapsto\phi^\vv_t$ on SDiff~$\Omega$. Suppose $t$ is fixed. After a small time $\tau$ the diffeomorphism describing the fluid will be $\exp_{Id}(\tau \vv(t))\phi_t^\vv$ up to second order terms in $\tau$. It implies $\vv(t)=DR_{(\phi_t^\vv)^{-1}}\dotex \phi_t^\vv$ where $DR_g$ denotes the tangent map induced by right multiplication by $g$ on the group. Thus the kinetic energy of the fluid  $T=\frac{1}{2}\langle\vv,\vv\rangle$ defines a right-invariant metric. The least action principle implies that the fluid motion  $t\mapsto\phi^\vv_t$ is a geodesic flow on SDiff~$\Omega$ endowed with the kinetic energy metric. Thus ${\nabla}^{LC}_\vv \vv=0$  where the Levi-Civita covariant differentiation ${\nabla}^{LC}$ is given by $\huge{\nabla}^{LC}_\vv\vec{\eta}=\frac{\partial}{\partial t}\xi+(\vv\cdot\nabla)~\vec{\eta}+\nabla{\alpha}$ and $\alpha$ is a real function such that  ${\nabla}^{LC}_\vv~\vec{\eta}\in\UU$. For  fixed $t$,
the virtual displacement corresponding to $\delta x$ in lemma 1 can be defined and identified to an element of $\UU$. It satisfies a Jacobi equation along the geodesic (see Proposition 2 of \cite{rouchon-EJM/B91}).

The reduced observer is defined intrinsically and can formally be applied to this fluid velocity estimation problem. The Theorem 1 is valid, as the proof is only made of intrinsic calculations, and its core is the Jacobi equation which gives conditions under which $\sigma \frac{\partial}{\partial \sigma}G\geq G$. The observer's state $\hat\xi$  is a virtual fluid, defined as a solution of  \eqref{obsreduit}, where $q$ is replaced by $\phi_t^\vv$.  Using the right group multiplication one can define  $\zeta(t)=\hat\xi(t)\circ (\phi_t^\vv)^{-1}$.  Note that $\zeta$ must remain in the group identity connected component so that \eqref{obsreduit} is well-defined.
\subsubsection{Discussion on the convergence and curvature}
When the curvature is bounded from above by $A=-B^2<0$, the geodesic
flow is sensitive to initial conditions, and admits ergodic
properties \cite{arnold-book-3}. Surprisingly, in this case the
observer is globally exponentially convergent by Theorem 1. When
there are always sections with negative curvature along a geodesic,
it is commonly assumed that the sensitivity to initial conditions is
still valid.

We have the following formal convergence result:  consider a
sinusoidal parallel stationary motion of a fluid in the tore
$T^2=\{(x,y),~x\mod2\pi,~y\mod2\pi\}$ given by the current function
$\psi=\cos(kx+ly)$ with $k,l\in\mathbb N$, and the velocity vector
field $\vec v=\text{rot}~\psi$.  Take $\hat{\xi}(0)={\phi^\vv_0}$
for   \eqref{obsreduit}. Then $\norm{\hat v(t)- \vec v}$ converges
exponentially to $0$. The proof is obvious as both points belong to
the same geodesic. But one can expect a great robustness to
measurement noise. Indeed \cite{arnold-book-3} proves the motion
defined by $\psi$ is a geodesic of SDiff $T^2$, and  the curvature
is non positive in \emph{all} planes containing $\vec v(x,y)$.
Moreover it is zero only in a family of planes of null measure. But
by Theorem 1 negative curvature implies global stability, and small
positive curvature implies a large basin of attraction.

More generally, Arnol'd  \cite{arnold-book-3} considers the group S$_0$Diff $T^2$ of diffeomorphisms preserving the center of gravity.  Calculations show  the curvature is positive ``only in a few sections". He suggests to consider the mean curvature along paths to characterize the stability of the flow. As a consequence, if the atmosphere was a bidimensional incompressible fluid on the earth viewed as $T^2$ (identify opposite sides of the planisphere),  the wind should be known up to 5 decimals  for a two-months weather's prediction.  Following this suggestion, as the curvature is positive in only in a few sections,  one could expect a good global behavior of the observer.

\section{Simulations on the sphere}Consider the inertial motion of a material point constrained to lie on the sphere $\SSS^2$. The speed is constant (see section \ref{holonomic:sec}) and assumed to be equal to $1$. One can always choose coordinates $q=(q_x,q_y,q_z)\in\RR^3$ such that the motion writes:
$
\dot q_x(t)=\cos (t),~
\dot q_y(t)=\sin (t),~
\dot q_z(t)=0
$.
Let $\hat \xi=(\hat \xi_x,\hat \xi_y,\hat \xi_z)\in\RR^3$. The observer equation \eqref{obsreduit} writes
\begin{align*}\dotex {\hat \xi_x}&= \frac{1}{\lambda}\varphi~\frac{((q\wedge\hat \xi)\wedge\hat \xi)_x}{\parallel (q\wedge\hat \xi)\wedge\hat \xi\parallel},\quad
\dotex {\hat \xi_y}=\frac{1}{\lambda}\varphi~\frac{((q\wedge\hat \xi)\wedge\hat \xi)_y}{\parallel (q\wedge\hat \xi)\wedge\hat \xi\parallel},\\
\dotex
{\hat \xi_z}&=\frac{1}{\lambda}\varphi~\frac{((q\wedge\hat \xi)\wedge\hat \xi)_z}{\parallel
(q\wedge\hat \xi)\wedge\hat \xi\parallel}
\end{align*}where $\lambda<0$ and $\varphi$ is the angle between $q$ and $\hat \xi$. As the geodesics of the sphere are great circles, $\varphi$ is the geodesic length between those two points. The inital conditions are : $q(0)=[1,0,0]^T$ and
$\hat \xi(0)=\frac{1}{\sqrt{2}}[0,1,1]^T$. To simulate the sensor's imperfections a white noise whose amplitude is  20$\%$ of the maximal value of the signal was added. $\hat \xi$ converges to the equator, and asymptotically follows $q$ at a distance $\abs{\lambda}$. The parallel transport $\hat v$ of $\dotex \hat \xi$ is an estimation of $v$ (not noisier than the measured signal). In fact for $\lambda<\pi/2$ the observer always converges in simulation, and for $\lambda>\pi/2$ it does not.

\begin{figure}
\includegraphics*[width=1\textwidth]{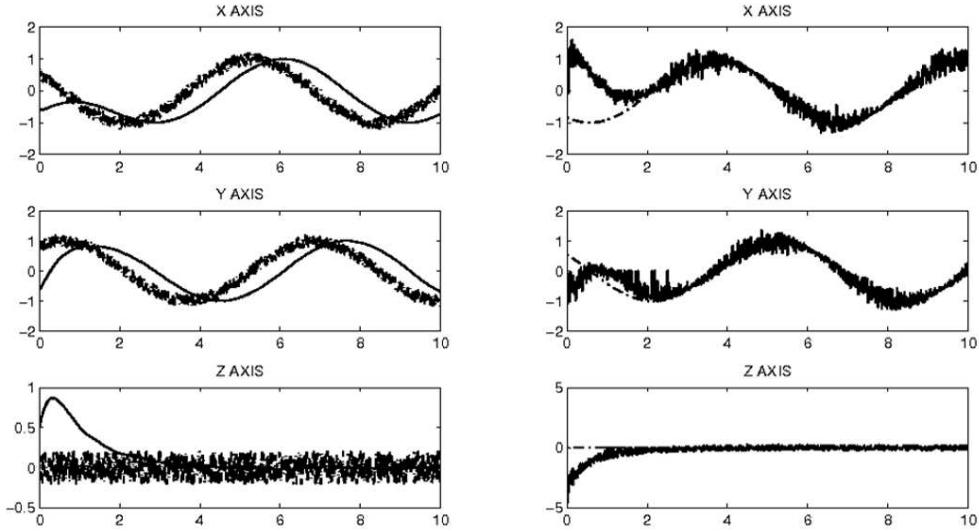}
  \caption{Simulations on the sphere for $\lambda=\pi/4$. Left: Measured $q$ (dashed line) and $\hat \xi$ (plain line). Right: velocity  $v$ (dashed) and estimation $\hat v$ (plain).}
    \label{simu:chariot1:fig}
\end{figure}
\section{Conclusion}We designed a nonlinear globally convergent reduced observer for conservative Lagrangian systems. The observer is intrinsic and converges despite the effects of curvature: instability of the flow and gyroscopic terms. The tuning of the  gains is simple. The only gain is a scalar which must be set in function of the noise and the maximal curvature. The observer can be used for velocity estimation for all systems described by geodesic flows ($\nabla_v v=0$), notably  conservative Lagrangian system, and the  motion of an incompressible fluid.   Using the Maupertuis principle this work could be extended to the case of a mixture of compressible fluids \cite{rouchon-ONERA1}.

Unfortunately when the motion is described by $\nabla_v v=S(q)$ with
$S$ known (Lagrangian system with external forces) the reduced
observer does not converge. Including such terms $S$ remains an open
question. As a concluding remark, note that the article gives
insight in the link between  convergence and geometrical structure
of the model in the theory of observers, complementing  the work of
\cite{aghannan-rouchon-ieee03,maithripala2005} and more recent results
\cite{arxiv-08}.
\section{Acknowledgements}
I am indebted to  Ivan Kupka and Pierre Rouchon  for very useful
discussions.

\bibliographystyle{plain}

\end{document}